\documentclass[reviewcopy]{elsart}

\usepackage{amssymb,amsmath}
\usepackage{graphicx}

\begin{document}

\begin{frontmatter}



\title{On an explicit finite difference method for fractional diffusion equations}


\author{S. B. Yuste\corauthref{sby}}
\corauth[sby]{Corresponding author.} \ead{santos@unex.es}
\ead[url]{http://www.unex.es/fisteor/santos/sby.html}
\author{ L. Acedo}

\address{Departamento de F\'{\i}sica, Universidad  de  Extremadura,
E-06071 Badajoz, Spain}

\begin{abstract}
A numerical method to solve the fractional diffusion equation,
which could also be easily extended to many other fractional
dynamics equations, is considered. These fractional equations have
been proposed in order to describe anomalous transport
characterized by non-Markovian kinetics and the breakdown of
Fick's law. In this paper we combine the  forward time centered
space (FTCS) method, well known for the numerical integration of
ordinary diffusion equations, with the Gr\"unwald-Letnikov
definition of the fractional derivative operator to obtain an
explicit fractional FTCS scheme for solving  the fractional
diffusion equation. The resulting method is amenable to a
stability analysis \emph{\`a la} von Neumann. We show that the
analytical stability bounds are in excellent agreement with
numerical tests. Comparison between exact analytical solutions and
numerical predictions are made.
\end{abstract}

\end{frontmatter}

\section{Introduction}
Fractional differential equations have been a highly specialized
and isolated field of mathematics for many years \cite{Podlubny}.
However, in the last decade there has been increasing interest in
the description of physical and chemical processes by means of
equations involving fractional derivatives and integrals. This
mathematical technique has a broad potential range of applications
\cite{HilferEd}: relaxation in polymer systems, dynamics of
protein molecules  and the diffusion of contaminants in complex
geological formations \cite{Kirchner,Scher,Berkowitz} are some of
the most recently suggested \cite{PhysToday}.

Fractional kinetic equations have proved particularly useful in
the context of anomalous slow diffusion (subdiffusion)
\cite{MetzlerRev}. Anomalous diffusion is characterized by an
asymptotic behavior of the mean square displacement of the form
\begin{equation}
\left< x^2(t)\right> \sim \frac{2K_\gamma}{\Gamma(1+\gamma)}
t^\gamma \; , \label{msd}
\end{equation}
where $\gamma$ is the anomalous diffusion exponent. The process is
usually referred to as subdiffusive when $0 < \gamma < 1$.
Ordinary (or Brownian) diffusion corresponds to $\gamma=1$ with
$K_1=D$ (the diffusion coefficient). From a continuous
(macroscopic) point of view, the diffusion process is described by
the diffusion equation $u_t(x,t)=D\, u_{xx}(x,t)$, where $u(x,t)$
represents the probability density of finding a particle at $x$ at
time $t$, and where $u_{\eta \zeta\ldots}$ is the partial
derivative with respect to the variables $\eta$,$\zeta$\ldots From
a microscopic point of view, the continuous description is known
to be connected with a Markov process in which the microscopic
particles (random walkers) perform stochastic jumps of finite mean
and finite variance. In these conditions the central limit theorem
holds for the sum of these jumps and Einstein's law for the mean
square displacement ensues [Eq.\ (\ref{msd}) with $\gamma=1$].

On the other hand, if an underlying non-Markovian microscopic
process is assumed in which random walkers perform jumps at times
chosen from a distribution with an algebraic long-time tail
$t^{-\gamma-1}$, then the diffusion process is anomalous
\cite{MetzlerRev,Rangarajan}.
In these circumstances the central limit theorem breaks down and
one must apply the generalized L\'evy-Gnedenko statistics
\cite{MetzlerRev,Hughes} which form the basis of Eq.\ (\ref{msd}).
It turns out that the probability density function $u(x,t)$ that
describes these anomalous diffusive particles follows the
fractional diffusion equation
\cite{MetzlerRev,VBPhysA,MellinTransform1,MellinTransform2}:
\begin{equation}
\label{subdeq} \frac{\partial }{\partial t} u(x,t)= K_\gamma
~_{0}D_{t}^{1-\gamma } \frac{\partial^2}{\partial x^2} u(x,t)
\end{equation}
where  $~_{0}D_{t}^{1-\gamma } $ is the fractional derivative
defined through the Riemann-Liouville operator (see Sec.\
\ref{sec:Def}). Fractional subdiffusion-advection equations, and
fractional Fokker-Planck equations  have also been proposed
\cite{MetzlerPRL,Advection1,Advection2,Advection3} and even
subdiffusion-limited reactions have been discussed within this
framework \cite{YusteKatja}.  In the mathematical literature,
these equations are usually referred to as parabolic
integro-differential equations with weakly singular kernels
\cite{Chen}.

These current applications of fractional differential equations
and many others that may well be devised in the near future make
it imperative to search for methods of solution. Some exact
analytical solutions for a few cases, although important, have
been obtained by means of the Mellin transform
\cite{MellinTransform1,MellinTransform2} and the method of images
\cite{images}. The powerful method of separation of variables can
also be applied to fractional equations in the same way as for the
usual diffusion equations (an example is given in Sec.\
\ref{sec:comp}). Another route to solving fractional equations is
through the integration of the product of the solution of the
corresponding non-fractional equation (the Brownian counterpart
obtained by setting $\gamma\rightarrow 1$) and a one-sided L\'evy
stable density \cite{MetzlerRev,subordination1,subordination2}.
However, as also for the Brownian case, the availability of
numerical methods for solving \eqref{subdeq} would be most
desirable, especially for those cases where no analytical solution
is available.
 One possibility was discussed recently by R. Gorenflo {\em et
al.} \cite{GorenfloA,GorenfloB,GorenfloC} who presented a scheme
to build discrete models of random walks suitable for the Monte
Carlo simulation of random variables with a probability density
governed by fractional diffusion equations.  Another more standard
approach is to build difference schemes of the type used for
solving Volterra type integro-differential equations \cite{Chen}.
In this line, some implicit (\emph{backward} Euler and
Crank-Nicholson) methods have been proposed
\cite{Chen,Serna,Lopez,LubichST,McLean1,McLean2,AdolfssonEL}. In
this paper we shall use the \emph{forward} Euler difference
formula for the time derivative $\partial u /\partial t$ in Eq.
\eqref{subdeq} to build an {\em explicit} method that we will call
the fractional Forward Time Centered Space (FTCS) method. For
Brownian ($\gamma=1$) diffusion equations, this explicit procedure
is the simplest numerical methods workhorse
\cite{MortonMayers,NumericalRecipes}. However, for fractional
diffusion equations, this explicit method has been overlooked
perhaps because of the difficulty in finding the conditions under
which the procedure is stable. This problem is solved here by
means of an analysis of Fourier--von Neumman type.

The plan of the paper is as follows. In Sec.\ \ref{sec:Def} we
give a short introduction to some  results and definitions in
fractional calculus. The numerical procedure to solve the
fractional diffusion equation \eqref{subdeq} by means of the
explicit FTCS method is given in Sec.\ \ref{sec:explicit}. In this
section we also discuss the stability and the truncating errors of
the FTCS scheme. In Sec.\ \ref{sec:comp} we compare exact
analytical solutions with the numerical ones and check the
reliability of the analytical stability condition. Some concluding
remarks are given in Sec.\ \ref{sec:Conclu}.

\section{Basic concepts of fractional calculus}
\label{sec:Def} The notion of fractional calculus was anticipated
by Leibniz, one of the founders of standard calculus, in a letter
written in 1695 \cite{Podlubny,MetzlerRev}. But it was in the next
two centuries that this subject fully developed into a field of
mathematics with work of Laplace, Cayley, Riemann, Liouville, and
many others.

There are two alternative definitions for the fractional
derivative $~_{0}D_{t}^{1-\gamma }$ of a function $f(t)$ which
coincide under relatively weak conditions. On the one hand, there
is the Riemann-Liouville operator definition
\begin{equation}
\label{RLdef} ~_{0}D_{t}^{1-\gamma } f(t)=\frac{1}{\Gamma(\gamma)}
\frac{\partial}{\partial t} \int_0^t d\tau
\frac{f(\tau)}{(t-\tau)^{1-\gamma}}\; ,
\end{equation}
with $0 < \gamma < 1$. For $\gamma=1$ one recovers the identity
operator and  for $\gamma=0$ the ordinary first-order derivative.
On the other hand, for any function $f(t)$ that can be expressed
in the form of a power series, the fractional derivative of order
$1-\gamma$ at any point inside the convergence region of the power
series can be written in the Gr\"unwald-Letnikov form
\begin{equation}
\label{GLdefpura} ~_{0}D_{t}^{1-\gamma }f(t)=\lim_{h \rightarrow
0} \frac{1}{h^{(1-\gamma)}} \sum_{k=0}^{[t/h]}
\omega_k^{(1-\gamma)} f(t-kh),
\end{equation}
where $[t/h]$ means the integer part of $t/h$. The
Gr\"unwald-Letnikov definition is simply a generalization of the
ordinary discretization formulas for integer order derivatives
\cite{Podlubny}. The Riemann-Liouville and the Gr\"unwald-Letnikov
approaches coincide under relatively weak conditions:  if $f(t)$
is continuous and $f'(t)$ is integrable in the interval $[0,t]$
then for every order $0 <1-\gamma < 1$ both the Riemann-Liouville
and the Gr\"unwald-Letnikov derivatives exist and coincide for any
time inside the interval $[0,t]$ \cite{Podlubny}. This theorem of
fractional calculus assures the consistency of both definitions
for most physical applications where the functions are expected to
be sufficiently smooth.

The Gr\"unwald-Letnikov  definition is important for our purposes
because it allows us to estimate $~_{0}D_{t}^{1-\gamma }f(t)$
numerically in a simple and efficient way:
 \begin{equation}
\label{GLdef} ~_{0}D_{t}^{1-\gamma }f(t)= \frac{1}{h^{(1-\gamma)}}
\sum_{k=0}^{[t/h]} \omega_k^{(1-\gamma)} f(t-kh) +O(h^p) \; ,
\end{equation}
The order of the resulting approximation, $p$, depends on the
choice of  $\omega_k^{(1-\gamma)}$. The approximation is of first
order ($p=1$) when $\omega_k^{(\alpha)}$ is the $k$-th coefficient
in the power series expansion of $(1-z)^\alpha$
\cite{Podlubny,LubichCoef}, i.e.,
\begin{equation}
\label{G1} (1-z)^\alpha=\sum_{k=0}^{\infty} \omega_k^{\alpha} z^k
\end{equation}
so that $\omega_k^{(\alpha)}=(-1)^k \binom{\alpha}{k}$ or,
equivalently:
\begin{equation}
\omega_0^{(\alpha)}=1, \qquad
\omega_k^{(\alpha)}=\left(1-\frac{\alpha+1}{k}
\right)\omega_{k-1}^{(\alpha)}\quad k=1,2,\ldots
 \label{coefO1}
\end{equation}
The coefficients $\omega_k^{(1-\gamma)}$ of the second-order
approximation ($p=2$) can be obtained similarly
\cite{Podlubny,LubichCoef}:
\begin{equation}
\label{G2} \left(\frac{3}{2}-2z+\frac{1}{2}
z^2\right)^\alpha=\sum_{k=0}^{\infty} \omega_k^{(\alpha)} z^k.
\end{equation}
These coefficients can be easily calculated using Fast Fourier
Transforms \cite{Podlubny}. However, for the fractional FTCS
method discussed in this paper, we will show in the next section
that nothing is gained by using second-order approximations for
the fractional derivative. Besides, the stability bound is smaller
if we take the coefficients derived from Eq.\ (\ref{G2}). Finally,
it is important to note that the error estimates given in
\eqref{GLdef} are valid only if either $t/h\gg 1$ \cite{Podlubny}
or $u(x,t)$ is sufficiently smooth at the time origin $t=0$
\cite{GorenfloNum}.

\section{Fractional Forward Time Centered Space method.}
\label{sec:explicit} We will use the customary notation $x_j=j
\Delta x$, $t_m=m\Delta t$ and $u(x_j,t_m)\equiv u_j^{(m)}\simeq
U_j^{(m)}$ where $U_j^{(m)}$ stands for the numerical estimate of
the exact value of $u(x,t)$ at the point $(x_j,t_m)$. In the usual
FCTS method, the diffusion equation is replaced by a difference
recurrence system for the quantities $u_j^{(m)}$:
\begin{equation}
\label{ec:difufw} \frac{u_j^{(m+1)}-u_j^{(m)}}{\Delta t}= D\,
\frac{u_{j-1}^{(m)}-2u_j^{(m)}+u_{j+1}^{(m)}}{(\Delta
x)^2}+T(x,t),
\end{equation}
with $T(x,t)$ being the truncation term \cite{MortonMayers}. In
the same way, the fractional equation is replaced by
\begin{equation}
\label{ec:met1} \frac{u_j^{(m+1)}-u_j^{(m)}}{\Delta t}= K_\gamma\,
~_{0}D_{t}^{1-\gamma}\frac{{u}_{j-1}^{(m)}-2{u}_j^{(m)}
+{u}_{j+1}^{(m)}}{(\Delta x)^2}+T(x,t) \; .
\end{equation}
The estimate of the truncation term will be given in Sec.\
\ref{subs:trunc}. Inserting the Gr\"unwald-Letnikov definition of
the fractional derivative given in Eq.\ (\ref{GLdef}) into Eq.\
(\ref{ec:met1}), neglecting the truncation term, and rearraging
the terms, we finally get the explicit FTCS difference scheme
\begin{equation}
\label{erme} U_j^{(m+1)}= U_j^{(m)}+S_\gamma \sum_{k=0}^{m}
\omega_k^{(1-\gamma)}
\left[U_{j-1}^{(m-k)}-2U_j^{(m-k)}+U_{j+1}^{(m-k)}\right] \; ,
\end{equation}
where $S_\gamma= K_\gamma {\Delta t}/[{h^{1-\gamma}(\Delta
x)^2}]$. In this scheme, $U_j^{(m+1)}$, for every position $j$, is
given \emph{explicitly} in terms of all the previous states
$U_j^{(n)}$, $n=0,1,\ldots,m$. Because the estimates $U_j^{(m)}$
of $u(x_j,t_m)$ are made at the times $m \Delta t$,
$m=1,2,\ldots$, and because the evaluation of
$~_{0}D_{t}^{1-\gamma} u(x_j,t)$ by means of \eqref{GLdef}
requires knowing $u(x_j,t)$  at the times $n h$, $n=0,
1,2,\ldots$, it is natural to choose $h=\Delta t$. In this case,
\begin{equation}
\label{sgam} S_\gamma= K_\gamma \frac{\Delta t^\gamma}{(\Delta
x)^2}\; .
\end{equation}
 The solution $u(x,t)$ is a causal function of
time with $u(x,t)=0$ if $t<0$ ($u_j^{(n)}=0$ if $n \le -1$), and
we assume that the system is prepared in an initial state
$u_j^{(0)}=U_j^{(0)}$. The iteration process described by Eq.\
(\ref{erme}) is  easily implementable as a computer algorithm, but
the resulting program is far more memory hungry than the
elementary Markov diffusive analogue because, in evaluating
$U_j^{(m+1)}$, one has to save all the previous estimates
$U_{j-1}^{(m+1)}$, $U_j^{(m+1)}$ and $U_{j+1}^{(m+1)}$ for
$n=0,1,\ldots m$. However, the use of the short-memory principle
\cite{Podlubny} could alleviate this burden. Anyway, before
tackling Eq.\ (\ref{erme}) seriously we must first discuss two
fundamental questions concerning any integration algorithm: its
stability and the magnitude of the errors committed by the
replacement of the continuous equation by the discrete algorithm.

\subsection{Stability of the fractional FTCS method}
\label{subs:stab} We will make a von Neumann type stability
analysis of the fractional FTCS difference scheme \eqref{erme}. We
start by assuming  a solution (a subdiffusion mode or
eigenfunction) with the form $u_j^{(m)}=\zeta_m e^{iqj\Delta x}$
where $q$ is a real spatial wave number. Inserting this expression
into \eqref{erme} one gets
\begin{equation}
\label{eqrigo} \zeta_{m+1}=\zeta_{m} -4 S
\sin^2\left(\frac{q\Delta x}{2}\right) \sum_{k=0}^{m}
\omega_k^{(1-\gamma)}\zeta_{m-k}\; .
\end{equation}
It is interesting to note that this equation is the discretized
version of
\begin{equation}
\frac{d \psi(t)}{d t}=- 4C \sin^2\left(\frac{q\Delta x}{2}\right)
~_{0}D_{t}^{1-\gamma}\psi(t) \; ,
\end{equation}
[with $C=S(\Delta t)^\gamma$] whose solution can be expressed in
terms of the Mittag-Leffler function $E_{\gamma}(-\lambda
t^{\gamma})$ \cite{HilferEd,MetzlerRev}. This result is not
unexpected because the subdiffusion modes of \eqref{subdeq} decay
as Mittag-Leffler functions \cite{MetzlerRev} [e.g., see
\eqref{uexact}].

\begin{figure}
\begin{center}
\includegraphics[width=0.95 \columnwidth]{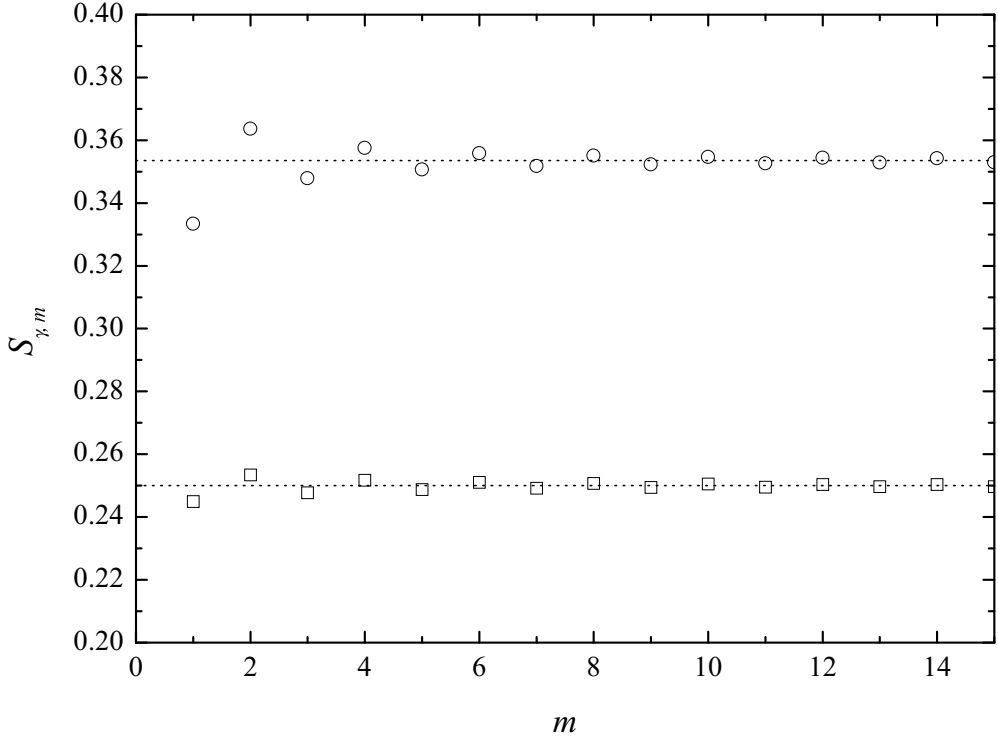}
\end{center}
\caption{First values of $S_{\gamma,m}$ versus $m$ for
$\gamma=1/2$ when the first-order coefficients (circles)  and
second-order coefficients (squares) are used. The lines mark the
corresponding limit values $S_{\gamma}^\times$ given by Eqs.\
\eqref{Sfirst} and \eqref{Ssecond} \label{fig1a}}
\end{figure}

\begin{figure}
\begin{center}
\includegraphics[width=0.95 \columnwidth]{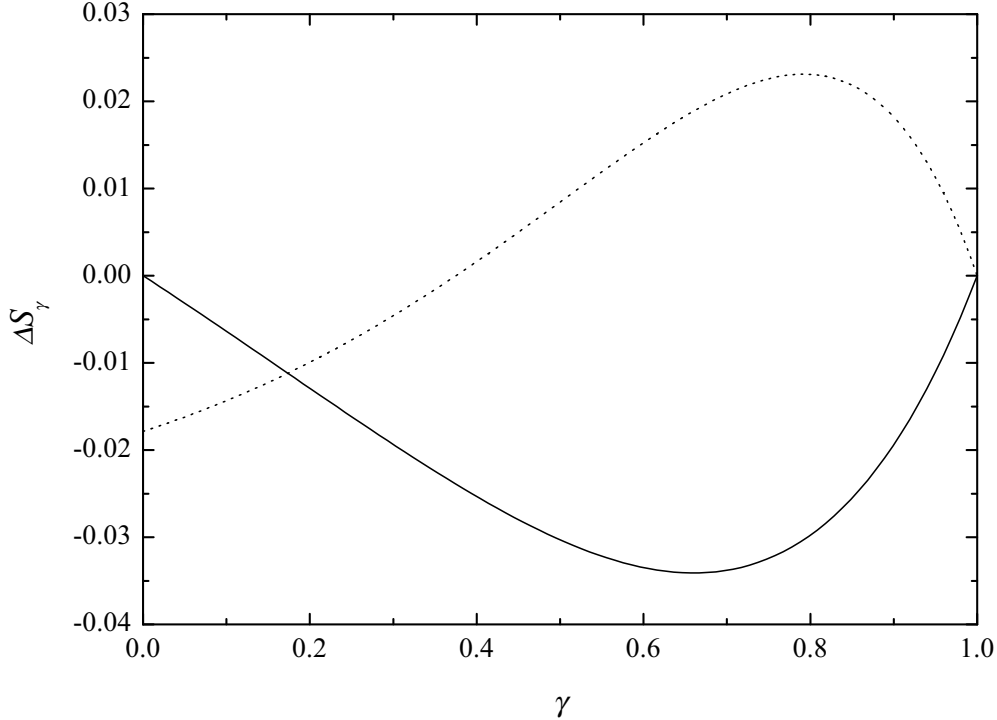}
\end{center}
\caption{The difference  $\Delta
S_{\gamma}=S_{\gamma,2}^{\times}-S_{\gamma,1}^{\times}$  versus
$\gamma$ when the first-order coefficients for
$\omega_k^{(1-\gamma)}$ [c.f. Eq. \protect{\eqref{coefO1}}] (solid
line) and second-order coefficients [c.f. Eq.
\protect{\eqref{G2}}] (dotted line) are used. \label{fig1b}}
\end{figure}

The stability of the solution is determined by the behaviour of
$\zeta_{m}$. Unfortunately, solving Eq.\ (\ref{eqrigo}) is much
more difficult than solvin the corresponding equation for the
diffusive case. However, let us write
\begin{equation}
\label{zxi} \zeta_{m+1}= \xi \zeta_m \; ,
\end{equation}
and let us \emph{assume } for the moment that $\xi\equiv\xi(q)$ is
independent of time.  Then  Eq.\ (\ref{eqrigo}) implies a closed
equation for the amplification factor $\xi$ of the subdiffusion
mode:
\begin{equation}\label{exr1}
\xi=1 -4 S_\gamma \sin^2\left(\frac{q\Delta
x}{2}\right)\sum_{k=0}^{m} \omega_k^{(1-\gamma)}\xi^{-k} \; .
\end{equation}
If $\vert \xi \vert > 1$ for some $q$, the temporal factor of the
solution grows to infinity according to Eq. (\ref{zxi}) and the
mode is unstable. Considering the extreme value $\xi=-1$, we
obtain from Eq.\ (\ref{exr1}) the following stability bound on
$S_\gamma$:
\begin{equation}\label{Sbound}
S_\gamma \sin^2\left(\frac{q\Delta x}{2}\right) \le
\frac{1/2}{\sum_{k=0}^{m}(-1)^{k} \omega_k^{(1-\gamma)}}\equiv
S_{\gamma,m}^{\times}\; .
\end{equation}
The bound expressed in Eq.\ (\ref{Sbound}) depends on the number
of iterations $m$. Nevertheless, this dependence is wak: for $m\ge
1$, $S_{\gamma,m}^{\times}$ approaches
$S_{\gamma,\infty}^{\times}\equiv S_{\gamma}^{\times}$ in the form
of   oscillations with small decaying amplitudes (see Fig.
\ref{fig1a}).  Figure \ref{fig1b}, in which we plot $\Delta
S=S_{\gamma,2}^{\times}-S_{\gamma,1}^{\times}$ versus $\gamma$ for
the first- and second-order coefficients, serves to gauge the
amplitude of these oscillations. In fact,  $\Delta S_\gamma$ is
the maximum value of
 $S_{\gamma,m+1}^{\times}-S_{\gamma,m}^{\times}$, $m\ge 1$ when the first-order
coefficients \eqref{coefO1} are used. We see that $\Delta
S_{\gamma}$ is certainly small for all $\gamma$.


The value of $\lim_{m\rightarrow \infty}S_{\gamma,m}^{\times} =
S_{\gamma}^{\times}$ can be deduced from Eq.\ (\ref{Sbound})
taking into account that the coefficients $\omega_k^{(1-\gamma)}$
are generated by the functions given in Eqs.\ (\ref{G1}) and
(\ref{G2}).  When the first-order coefficients given by \eqref{G1}
are used, one gets:
\begin{equation}
\label{Sfirst} S_\gamma^\times=
\frac{1}{\left.2(1-\xi)^{1-\gamma}\right|_{\xi\rightarrow
-1}}=\frac{1}{2^{2-\gamma}}\; .
\end{equation}
Similarly, when the second-order coefficients given by \eqref{G2}
are used, one gets:
\begin{equation}
\label{Ssecond} S_\gamma^\times=\frac{1}
{\left.
2\left(\frac{3}{2}-2\xi+\frac{1}{2}
\xi^2\right)^{1-\gamma}\right|_{\xi\rightarrow -1}}
=\frac{1}{4^{3/2-\gamma}}.
\end{equation}
We will verify numerically in Sec.\ \ref{sec:comp} that the
explicit integration method as given by Eq.\ (\ref{erme}) is
stable when
\begin{equation}\label{}
 S_\gamma \le \frac{S_\gamma^\times}{\sin^2\left(\frac{q\Delta
x}{2}\right)}
\end{equation}
and unstable otherwise.   As the maximum value of the square of
the sine function is bounded by 1, we can give a more conservative
but simpler bound: the fractional FTCS method will be stable when
\begin{equation}\label{SboundSimple}
 S_\gamma = K_\gamma \frac{\Delta t^\gamma}{(\Delta
x)^2}\le {S_\gamma^\times}.
\end{equation}
The physical interpretation of this restriction is the same as for
the diffusive case, namely, Eq.\ \eqref{SboundSimple} means that
the maximum allowed time step $\Delta t$ is, up to a numerical
factor, the (sub)diffusion time across a distance  of length
$\Delta x$  [c.f. Eq.\ \eqref{msd}].

Notice that the value of $S_\gamma^\times=1/4^{3/2-\gamma}$ given
by Eq.\ (\ref{Ssecond}) is smaller than $1/2^{2-\gamma}$ for any
$\gamma < 1$ (if $\gamma=1$ we recover the bound $S^\times=1/2$ of
the usual explicit FTCS method for the ordinary diffusion equation
\cite{MortonMayers,NumericalRecipes}). Consequently, the
fractional FTCS method that uses a second-order approximation in
the fractional derivative is ``less robust'' than the fractional
FTCS method that uses the first-order coefficients
$\omega_k^{(1-\gamma)}$. Taking into account that the two methods
have the same precision (see Sec.\ \ref{subs:trunc}) we note that
nothing is gained by using the fractional derivative with higher
precision. Therefore, in practical applications, we will only use
here the first-order coefficients \eqref{coefO1}.

\subsection{Truncating error of the fractional FTCS method}
\label{subs:trunc}

The truncating error $T(x,t)$ of the  fractional FTCS difference
scheme is [see \eqref{ec:met1}]:
\begin{equation}
\label{txtE} T(x,t)=\frac{u_j^{(m+1)}-u_j^{(m)}}{\Delta t}-
K_\gamma\,D_{t}^{1-\gamma}
\left[\frac{{u}_{j-1}^{(m)}-2{u}_j^{(m)} +{u}_{j+1}^{(m)}}{(\Delta
x)^2}\right] \; .
\end{equation}
But
\begin{equation}\label{opdifddt}
\frac{u_j^{(m+1)}-u_j^{(m)}}{\Delta t}= u_t+\frac{1}{2} u_{tt}
\Delta t+O(\Delta t)^2
\end{equation}
and
\begin{equation}\label{opdiffrac}
~_{0}D_{t}^{1-\gamma }\left[{{u}_{j-1}^{(m)}-2{u}_j^{(m)}
+{u}_{j+1}^{(m)}}\right]= \displaystyle\frac{1}{h^{1-\gamma}}
\sum_{k=0}^m w_k^{1-\gamma} \left[u_{xx}+\displaystyle\frac{1}{12}
u_{xxxx} \left(\Delta x \right)^2+\cdots\right]
+{O}\left(h^p\right)\;
\end{equation}
so that, taking into account that $u(x,t)$ is the exact solution
of Eq.\ (\ref{subdeq}), we finally get from Eqs.\ (\ref{txtE}),
(\ref{opdifddt}) and (\ref{opdiffrac}) the following result
\begin{align}
\label{txtF1} T(x,t)&={O}(h^p)+\frac{1}{2}u_{tt}\Delta t
-\frac{K_\gamma(\Delta x)^2}{12}~_{0}D_{t}^{1-\gamma }
u_{xxxx}+\cdots \\
&={O}(h^p)+ {O}(\Delta t)+ {O}(\Delta x)^2 \;. \label{txtF2}
\end{align}

Therefore, (i) assuming that the initial boundary data for $u$ are
consistent (as assumed for the usual FTCS method
\cite{MortonMayers}) and (ii) assuming that  $u$ is sufficiently
smooth at the origin $t=0$ [see remark below Eq.\ \eqref{G2}], we
conclude that the method discussed in this paper is
unconditionally consistent for any order $p$ because
$T(x,t)\rightarrow 0 $ as $h$, $\Delta t$, $\Delta x\rightarrow
0$.  As remarked above, in practical calculations is convenient to
use $h=\Delta t$ so that, due to the term ${O}(\Delta t)$ in
\eqref{txtF2}, no improvements are achieved by considering higher
orders than $p=1$ in the fractional derivative. In is interesting
to note that for the diffusion equation ($\gamma=1$) it is
possible to cancel out the last two terms in Eq.\ (\ref{txtF1})
with the choice $\Delta t=(\Delta x)^2/(6 K_\gamma)$, trhereby
obtaining a scheme that is ``second-order accurate''
\cite{MortonMayers}. This is not possible for the fractional case
because of the fractional operator.

\section{Numerical solutions and the stability bound on $S_\gamma$}
\label{sec:comp} The objective of this section is twofold: first
we want to test the reliability of the numerical algorithm defined
in Eq.\ (\ref{erme}) by applying it to two fractional problems
with known exact solutions, and second we want to check the
stability bounds obtained in Sec.\ \ref{subs:stab}.

\begin{figure}
\begin{center}
\includegraphics[width=0.95 \columnwidth]{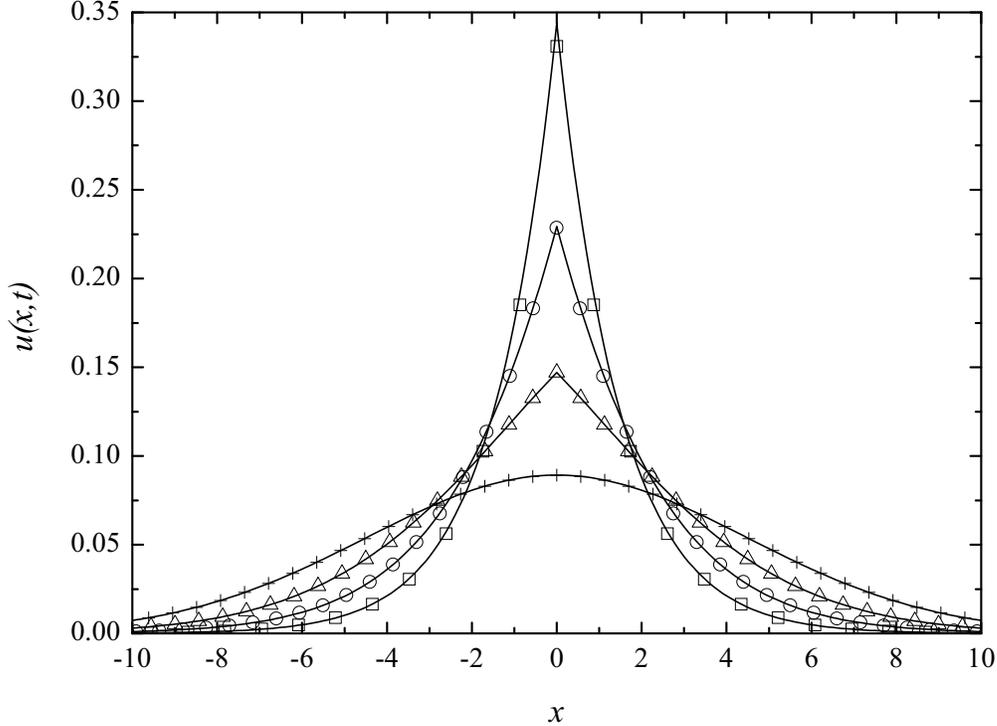}
\end{center}
\caption{Comparison between the exact subdiffusion propagator
(lines) and the numerical integration results for $\gamma=1/4$
(squares), $\gamma=1/2$ (circles), $\gamma=3/4$ (triangles) and
$\gamma=1$ (crosses) and $t=10$.
\label{fig1}}
\end{figure}

\subsection{Numerical solution versus exact solution: two examples}
The  fundamental solution of the subdiffusion equation in Eq.\
(\ref{subdeq}) corresponds to the problem defined in the unbounded
space where the initial condition is $u(x,t=0)=\delta(x)$. This
solution is called the propagator (or Green's function) and can be
expressed in terms of Fox's H-function \cite{MetzlerRev}:
\begin{equation}
\label{propag} u(x,t)=\frac{1}{\sqrt{4\pi K_\gamma t^\gamma}}\,
H^{10}_{11}\left[\frac{|x|}{\sqrt{ K_\gamma t^\gamma}}
\left|\begin{array}{l}{(1-\gamma/2,\gamma/2)}\\[1ex]{(0,1)}\end{array} \right.
\right]\; .
\end{equation}
In our numerical solution we used the boundary conditions
$u(-L,t)=u(L,t)=0$ with a sufficiently large $L$ in order to avoid
finite size effects. In Fig.\ \ref{fig1} we compare the numerical
integration results with the exact solution (\ref{propag}) for
 $\gamma=1/4$, $1/2$, $3/4$, $1$ at
$t=10$. The timestep used was $\Delta t=0.01$ and $\Delta
x=\sqrt{K_\gamma (\Delta t)^\gamma/S_\gamma}$ with $K_\gamma=1$
and $S_\gamma=0.28$, $0.33$, $0.4$ and $0.5$. All these values of
$S_\gamma$  are just below  the stability bound $S_\gamma^\times$
(see Eq.\ (\ref{Sfirst})). The agreement is excellent except for
$\gamma=1/4$ and $x=0$, but this minor discrepancy  is surely due
 to the large spatial cell $\Delta x\simeq 1.06$ used in this case.

\begin{figure}
\begin{center}
\includegraphics[width=0.95 \columnwidth]{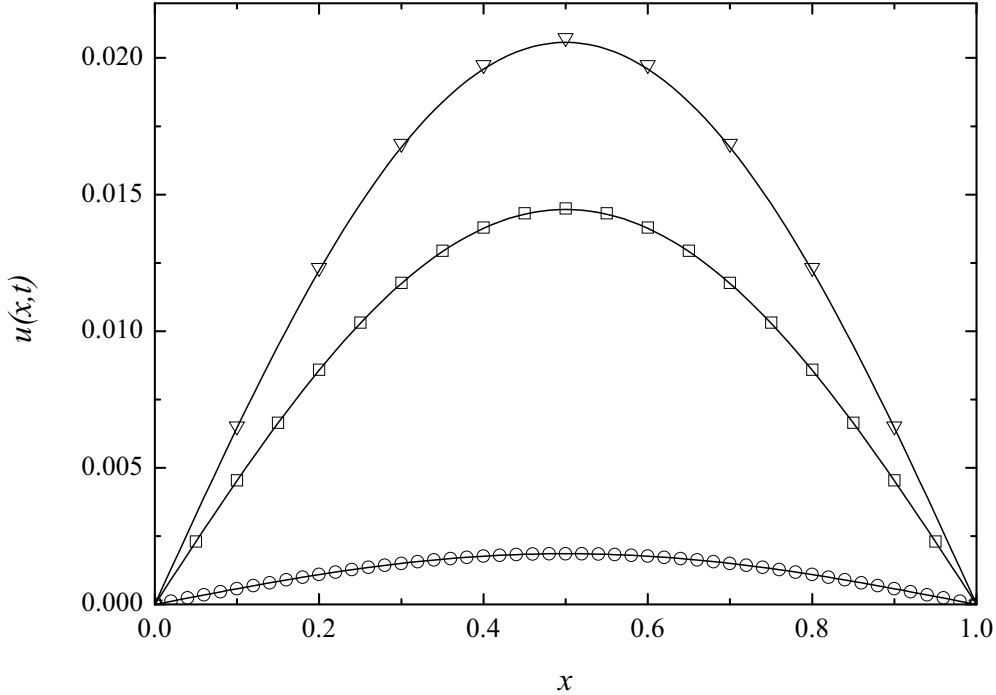}
\end{center}
\caption{Numerical solution of the subdiffusion equation for the
problem with absorbing boundary conditions, $u(0,t)=u(1,t)=0$, and
initial condition $u(x,0)=x(1-x)$ versus the exact analytical
result (lines) for $t=0.5$. The solution $u(x,t)$ is shown  for
$\gamma=0.5$ (triangles), $\gamma=0.75$ (squares) and $\gamma=1$
(circles). \label{fig2}}
\end{figure}
We have also considered a problem with absorbing boundaries,
$u(0,t)=u(1,t)=0$, and initial condition $u(x,t=0)=x(1-x)$. The
exact analytical solution of Eq.\ (\ref{subdeq}) is easily found
by the method of separation of variables: $u(x,t)=X(x)T(t)$. We
thus find $X_n(x)=\sin(n \pi x)$ and
\begin{equation}
\label{Tsepeq} \frac{dT}{dt}=-K_\gamma\, \lambda_n^2
~_{0}D_{t}^{1-\gamma } T\; ,
\end{equation}
where $\lambda_n=n \pi$, $n=1,2,\ldots$. The solution of Eq.\
(\ref{Tsepeq}) is found in terms of the Mittag-Leffler function
\cite{MetzlerRev}:
\begin{equation}
\label{Tn} T_n(t)=E_\gamma(-K_\gamma n^2 \pi^2 t^\gamma)\; .
\end{equation}
Imposing the initial condition we obtain
\begin{equation}
\label{uexact} u(x,t)=\frac{8}{\pi^3} \sum_{n=0}^\infty
\frac{1}{(2n+1)^3} \sin[(2n+1)\pi x] E_\gamma[-K (2n+1)^2\pi^2
t^\gamma]\; .
\end{equation}
In Fig.\ \ref{fig2} we compare this exact solution with the
results of the numerical integration scheme for $\gamma=0.5$,
$\gamma=0.75$, and $\gamma=1$ for $t=0.5$ and $K_\gamma=1$. The
values of $S_\gamma$ used were $S_\gamma=0.33$, $0.4$, and $0.5$
with $\Delta x=1/10$, $1/20$, and $1/50$, respectively. The values
of $\Delta t$ for fixed $S_\gamma$ and $\Delta x$ stem from the
definition of $S_\gamma$:
\begin{equation}
\label{Sdef} \Delta t=\left[ \frac{S_\gamma(\Delta x)^2}{K_\gamma}
\right]^{1/\gamma}\; .
\end{equation}
Excellent agreement is observed for the three values of $\gamma$,
it being slightly poorer for the smallest value which is not
surprising because in this case the mesh size $\Delta x=1/10$ used
is the largest.

\begin{figure}[t]
\begin{center}
\includegraphics[width=0.95 \columnwidth]{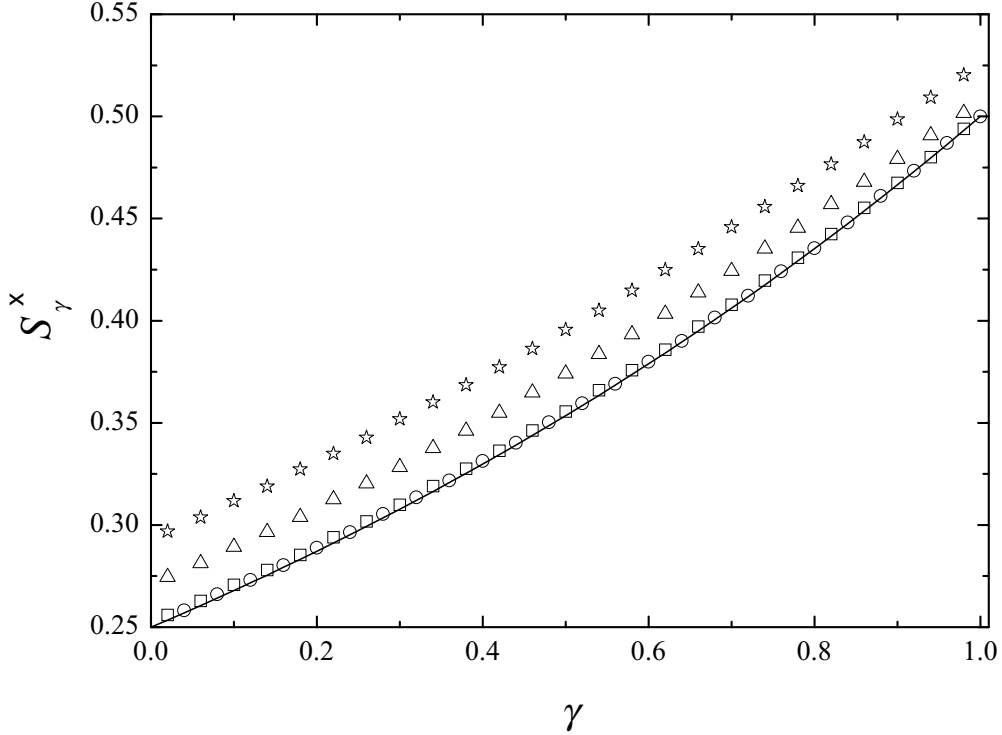}
\caption{Values of $S_\gamma^\times$ corresponding to the onset of
instability versus the subdiffusion exponent $\gamma$. The solid
line is the prediction of the Fourier--von Neumann analysis and
the symbols denote the results of the numerical tests with the
criterion in Eq.\ (\protect\ref{scriter}): stars, triangles and
squares for the absorbing boundary problem with $u(x,0)=x(1-x)$
with $M=50$, $100$ and $1000$, respectively, and circles for the
propagator with $M=1000$.
\label{fig3}}
\end{center}
\end{figure}

\subsection{Numerical check of the stability analysis}
We  checked  the stability bound on the value of the $S_\gamma$
given in Eq.\ (\ref{Sfirst}) in the following way.  For a set of
values of $\gamma$ in the interval $[0,1]$, and for values of
$S_\gamma$ starting at $S_\gamma=0.98 S_\gamma^\times$ (in
particular, for $S_\gamma=0.98/ 2^{2-\gamma}+ 0.001\,n$,
$n=0,1,2,\ldots$) we  applied the fractional FTCS  integration
until step $M$. We say that the resulting integration for a given
values of $\gamma$ and $S_\gamma$ is unstable when the following
condition is satisfied at any position $j$:
\begin{equation}
\label{scriter} \left\vert\frac{u_j^{m-1}}{u_j^m}-\Xi\right\vert >
\Xi \quad \text{for any} \quad m=M-\Delta M,M-\Delta
M+1,\ldots,M\; ,
\end{equation}
where $\Xi=5$ and $\Delta M=10$. This means that the numerical
solution is considered unstable if the quotient
$u_j^{m-1}/u_j^{m}$ becomes negative or larger than $2 \Xi$ at any
of the last $\Delta M$ steps. (Of course, this criterion is
arbitrary; however, the results do not change substantially for
any other reasonable choice of  $\Xi$ and $\Delta M$.) Let
$S_\gamma^\text{min}$ be the smallest value of $S_\gamma=0.98/
2^{2-\gamma}+ 0.001\,n$ that verifies the criterion
(\ref{scriter}). For the absorbing boundary problem  we calculate
these values using $\Delta x=1/2N$ with $N=5$ and $M=50$, $M=100$
and $M=1000$. For the propagator, we calculate
$S_\gamma^\text{min}$ using $M=1000$ and $\Delta t=5\times
10^{-4}$ in a lattice with absorbing frontiers placed at
$x=-N\Delta x$ and $x=N\Delta x$ with $N=50$. It is well known
that for a lattice with $2N+1$ points (including the absorbing
boundaries) the maximum value of $\sin(q\Delta x/2)$ in Eq.\
\eqref{Sbound} occurs for $q\Delta x=(2N-1)\pi/(2N)$, so that in
Fig.\ (\ref{fig3}) we  plot  $S_\gamma^\text{min}
\sin^2[(2N-1)\pi/(4N)]$. We observe that for large $M$ the
stability bound predicted by  Eq.\ (\ref{Sfirst}) agrees with the
result of the numerical test. The larger values obtained for
smaller $M$ mean that the method must be ``very unstable'' to
fulfill our instability criterion in so few steps. The success of
the numerical test is truly remarkable and supports the unorthodox
application of the Fourier--von Neumann stability analysis to the
fractional FTCS scheme made in Sec.\ \ref{subs:stab}.

\begin{figure}
\begin{center}
\includegraphics[width=0.95 \columnwidth]{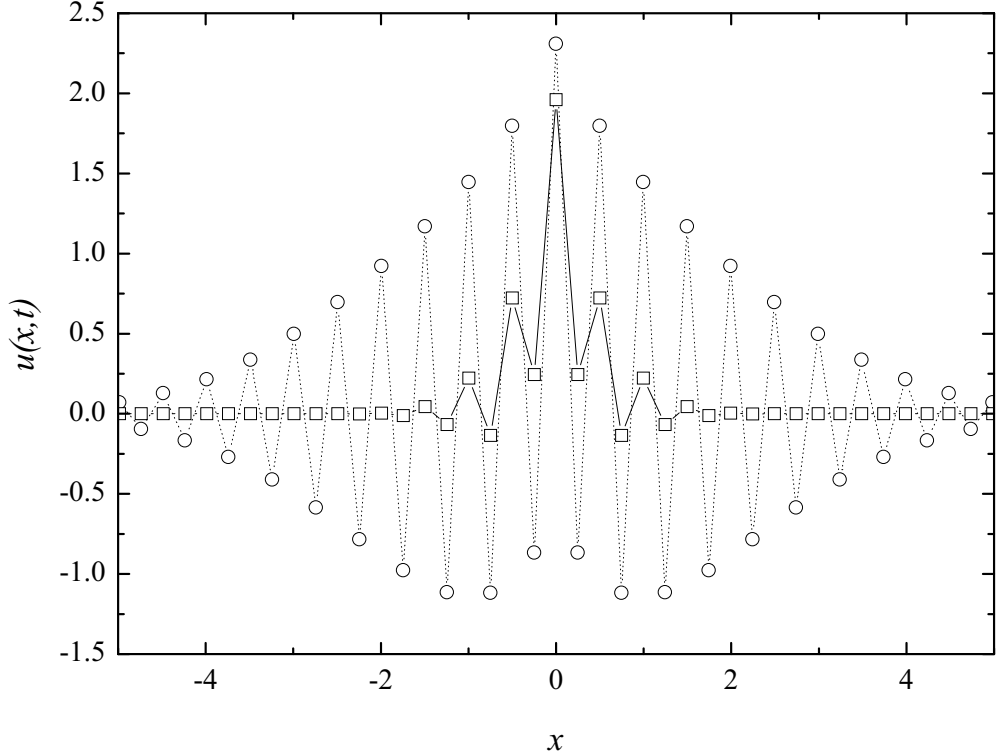}
\end{center}
\caption{The propagator $u(x,t)$ for $\gamma=1/2$, $K_\gamma=1$,
$S=0.36$ and $t=0.005$ (squares) and $t=0.05$ (circles). The time
step is $\Delta t=0.0005$ and the spatial mesh $\Delta x$ is
obtained according to Eq.\ (\protect\ref{Sdef}). The lines are
plotted as a visual guide. \label{fig4}}
\end{figure}

In Fig.\ (\ref{fig4}) we plot the numerical solution when
$S_\gamma=0.36 > S_\gamma^\times$ in the case of the propagator
with $\gamma=1/2$. This kind of oscillatory behaviour in the
unstable domain is typical for ordinary partial differential
equations too.

\section{Concluding remarks}
\label{sec:Conclu}

The availability of efficient numerical algorithms for the
integration of fractional equations is important as these
equations are becoming essential tools for the description of a
wide range of systems \cite{PhysToday}. In this paper we have
discussed a numerical algorithm for the solution of the fractional
(sub)diffusion equation \eqref{subdeq}. Although  we have dealt
with this particular equation, our procedure could be extended to
any fractional integro-differential equation by means of an
obvious combination of the Gr\"unwald-Letnikov definition of the
fractional derivative \cite{Podlubny,HilferEd,MetzlerRev} with
standard discretization algorithms used in the context of ordinary
partial differential equations \cite{MortonMayers}. Furthermore,
the method (given its explicit nature) can be trivially extended
to $d$-dimensional problems, which is not such an easy task when
implicit methods are considered.

In our numerical method the state of the system at a given time
$t=m \Delta t$ is given explicitly in terms of the previous states
at $t=(m-1) \Delta t, \ldots,\Delta t, 0$ by means of the FTCS
scheme \eqref{erme}. We verified that for some standard initial
conditions with exact analytical solution, namely, (a) the
propagator in an unlimited system with $u(x,t=0)=\delta(x)$ and
(b) a system with absorbing boundaries and $u(x,t=0)=x(1-x)$, the
present algorithm leads to numerical solutions which are in
excellent agreement with the exact solutions. Using a Fourier--von
Neumann technique we have provided the conditions for which the
fractional FTCS method is stable. For example, if a first-order
approximation for the fractional derivative is considered, we have
shown that the FTCS algorithm is stable if $S_\gamma=K_\gamma
(\Delta t)^\gamma/(\Delta x)^2\leq 1/2^{2-\gamma}$. For $\gamma=1$
the well-known bound $S=D \Delta t/(\Delta x)^2 \leq 1/2$ of the
ordinary explicit method for the diffusion equation is recovered.

Concerning the implementation of the method we must remark that
the evaluation of the state of the system at a given time step $m
\Delta t$ requires information about all previous states at
$t=(m-1) \Delta t, (m-2)\Delta t,\ldots,\Delta t, 0$ and not
merely the immediately preceding one as occurs in ordinary
diffusion. This is a consequence of the non-Markovian nature of
subdiffusion and implies the need for massive  computer memory in
order to store the evolution of the system, which is especially
cumbersome in computations of long-time asymptotic behaviours.
This could be palliated by using the ``short-memory'' principle
\cite{Podlubny}. Another feature of the explicit numerical scheme
is the interdependence of the temporal and spatial discrete steps
for a fixed $S_\gamma$. If, as usual, one intends to integrate an
equation with a given mesh $\Delta x$, then the corresponding step
size $\Delta t$  for a given  $S_\gamma < S_\gamma^\times$ is of
the order $(\Delta x)^{2/\gamma}$. As a consequence, $\Delta t$
could become extremely small even for no too small values of
$\Delta x$, especially when the problem is far from the diffusion
limit, i.e., for small values of $\gamma$, so that the number of
steps needed to reach even moderate times would become
prohibitively large. In this case,  the resort to implicit methods
\cite{Chen,Serna,Lopez,LubichST,McLean1,McLean2,AdolfssonEL},
stable for any value of $\Delta t$ and $\Delta x$, is compulsory.

This work has been supported by the Ministerio de
Ciencia y Tecnolog\'{\i}a (Spain) through Grant No. BFM2001-0718.



\end{document}